\documentclass[12pt]{article}
\usepackage{amssymb}
\usepackage{amsfonts}
\usepackage{latexsym}
\def\appendix#1{
\addtocounter{section}{1}
\renewcommand{\thesection}{\Alph{section}}
\section*{Appendix \thesection\protect\indent #1}
}
\def\egyenk{\begin{eqnarray}}
\def\egyenv{\end{eqnarray}}
\def\ds{\displaystyle}
\author{Zolt\'an Nagy \and Jean Avan \and Genevi\`eve Rollet}
\title{Construction of dynamical quadratic algebras}
\newtheorem{L1}{Lemma}
\newtheorem{T1}{Theorem}
\newtheorem{P1}{Proposition}
\begin{document}

\begin{titlepage}
%\setcounter{page}{1}
%\renewcommand{\thefootnote}{\fnsymbol{footnote}}

%\begin{flushright}

%\end{flushright}

\vspace{5 mm}

\begin{center}
{\Large\bf Construction of dynamical quadratic algebras}

\vspace{10 mm}

{\bf Zolt\'an Nagy\footnote{e-mail: nagy@ptm.u-cergy.fr}},
{\bf Jean Avan\footnote{e-mail: avan@ptm.u-cergy.fr}}
\textrm{and}
{\bf Genevi\`eve Rollet\footnote{e-mail: rollet@ptm.u-cergy.fr }}

\vspace{15mm}

{\it Laboratory of Theoretical Physics and Modelization\\
     University of Cergy-Pontoise (CNRS UMR 8089), 5 mail Gay-Lussac, Neuville-sur-Oise,\\
     F-95031 Cergy-Pontoise Cedex, France}
\end{center}

\vspace{18mm}

\begin{abstract}
\noindent We propose a dynamical extension of the quantum 
quadratic exchange 
algebras introduced by Freidel and Maillet. It admits two distinct fusion structures.
A simple example is provided by the scalar Ruijsenaars-Schneider model.
\end{abstract}

\vfill

\end{titlepage}

\section{Introduction}

The notion of dynamical quantum algebra was introduced by Felder \cite{Fe} and Gervais-Neveu \cite{GN} and subsequently 
studied and developed in
\cite{Etingof,Schiffmann,Enriquez}. The classical limit (dynamical r-matrices) first appeared in \cite{AT} 
and was later investigated in \cite{Schiffmann2}; examples related to the Ruijsenaars-Schneider (RS) model \cite{RS} 
were particularly studied \cite{AR, Suris}. In general, they are characterized by the existence, in the quantum
(and classical) $R$-matrices, of supplementary parameters identified as coordinates on the dual of
some particular Lie algebra $\mathfrak{h}$. These parameters occur as dynamical variables in the classical case,
hence the name.

These investigations until now concentrated on particular dynamical extensions of the quantum group structure
characterised by its quadratic exchange relation:
\egyenk 
R_{12} \ T_1 \ T_2 = T_2 \ T_1 \ R_{12}\nonumber
\egyenv
They were recently understood as Drinfel'd twists of the quantum group \cite{japonais}. The quantum $R$-matrix 
obeys a dynamical cubic equation (Gervais-Neveu-Felder (GNF) equation), generalizing \cite{Fe, GN} 
the quantum Yang-Baxter (YB) equation \cite{McGu,Ber,BZinn,Ya,Ba}
%\cite{GN,Fe}.

We will describe here an extension to quadratic exchange algebras 
of the notion of quantum dynamical algebra. 

These algebras are characterised by so-called braided exchange relations \cite{MFr}:
\egyenk
A_{12} \ T_1 \ B_{12}\ T_2 = T_2\ C_{12}\ T_1\ D_{12}\label{QQA}
\egyenv
where the generators of the algebra sit in the entries of $T$ viewed as a matrix in
$\mathrm{End}(V)$ for a given ``auxiliary'' vector space $V$; $A$, $B$, $C$, $D$ are c-number
structure matrices acting on $V \otimes V$. V may have the structure of a loop space $V \otimes C(\lambda)$ in which case the structure matrices depend on complex spectral parameters $\lambda_1, \lambda_2$. As usual in this context, the indices $1,2$ label the auxiliary vector spaces 
that the matrices act on.
%In the sequel, we refer to the collection of these four matrices as {\it the} quantum quadratic $R$-matrix.
Many examples are known in the case where  $A, B, C, D$ depend only on spectral parameters, see \cite{BA}.
Recently, a universal structure was proposed for the specific case of reflection algebras $A=C=B^{\pi}=D^{\pi}$ \cite{DKM}. 
Using associativity to compare both ways of exchanging $T_1 T_2 T_3$ into $T_3 T_2 T_1$ leads naturally
to YB-type equations on  the $A,B,C,D$ matrices as sufficient consistency conditions. 
In \cite{MFr} Maillet and Freidel wrote down 8 YB-equations 
which provide a sufficient set of 3-space exchange
conditions for the quantum algebra (\ref{QQA}). This case is hereafter refered to as `` nondynamical''.

The question now arises whether there exists a consistent way of 
dynamizing these 8 YB-equations (in the sense of Gervais-Neveu-Felder \cite{Fe, GN}) and whether such dynamized 
quantum YB-equations can be interpreted
as sufficient 3-space exchange conditions for a dynamical quadratic quantum algebra.

For simplicity we choose $\mathfrak{gl}(n)$ as underlying Lie algebra and its Cartan subalgebra as the Lie algebra
supporting the dynamical parameters.
We will here define a dynamical quadratic quantum algebra (DQQA) for a particular choice  of zero-weight conditions 
of the R-matrix set $A,B,C,D$ under the action of the Cartan subalgebra $\mathfrak{h}$. 
This choice is consistent --
as we will see -- with the specific structure of the classical and quantum $R$-matrices for RS-models \cite{Suris} 
and provides the general algebraic frame for
the construction of quantum RS-models proposed by Arutyunov-Chekov-Frolov (ACF) \cite{russes}. 
We will prove that such sufficient conditions for 3-space exchange of these DQQA's realize 
a dynamical version of the quadratic YB-equations
derived  in \cite{MFr}. We will also describe the ACF example of DQQA coming from the scalar RS model.
We will then describe two distinct coproduct-type structures for this  DQQA, generalizing 
the coproduct structures described in \cite{MFr}; these coproducts allow for building other spin-chain type models
from the scalar one.
We will finally define a classical limit of the DQQA and show that
the scalar Ruijsenaars-Schneider classical $r$-matrix structure does realize this classical limit (see \cite{russes}).
% 
%, thereby providing us with the first
% fully algebraic description of the classical YB-equations for this $r$-matrix and for the scalar Calogero-Moser (CM) $r$-matrix.

\section{Dynamical quadratic quantum algebras}

We start by expliciting the ``dynamical'' notation.
Let $\mathfrak{g}$ be a simple Lie algebra and $\mathfrak{h}$  a  commutative subalgebra of $\mathfrak{g}$ of dimension $n$. 
(For an extension to noncommutative  $\mathfrak{h}$ see \cite{Ping}.)
%f $\mathfrak{h}$ by $x$. 
Let us choose a basis $\ds\{h^{i}\}_{i=1}^n$ of $\mathfrak{h}^{\ast}$ and 
%let $\{\lambda^i\}_{i=1}^n$
let $\ds \lambda =\sum_i \lambda_i h^i$ be an element of $\mathfrak{h}^{\ast}$. 
The dual basis is denoted in $\mathfrak{h}$ by $\ds\{h_{i}\}_{i=1}^n$.
For any differentiable function
$f(\lambda)=f(\{\lambda_i\})$ one defines:
%taking values in $\mathbb{C}$, we denote
\egyenk
f(\lambda+\gamma h)= e^{\gamma \mathcal{D}}f(\lambda)e^{-\gamma \mathcal{D}},
\egyenv 
where
\egyenk
\mathcal{D}=\sum_i h_i \partial_{\lambda_i}
\egyenv
%be the coordinate system associated with this base in $\mathfrak{h^{\ast}}$. 
%For a matrix-valued function $R : \mathfrak{h^{\ast}} \to \textrm{End} (V \otimes V)$, where $V$ is a finite dimensional
%vector space that is also an $\mathfrak{h}$-module, $R_{12}(\lambda+\gamma h_3)$ is defined by the formula
%$R_{12}(\lambda+\gamma h_3)(v_1 \otimes v_2 \otimes v_3) = \big(R_{12}(\lambda +\gamma \mu) \otimes \mathbf{1}\big)
%(v_1 \otimes v_2 \otimes v_3)$ if $v_3$ is of weight $\mu$.
%
It can be seen that this definition yields formally 
%
%Finally, let us suppose that every vector space 
%appearing later, e.g. $U_1, U_2$, is also  
%let $V_1, \dots , V_m$ be 
%a finite-dimensional semisimple 
%$\mathfrak{h}$-module.
%Dynamical algebras are defined by relations similar to (\ref{QQA}) but  with  matrices depending on a parameter
%in a certain $\mathfrak{h}$. Thus, YB-equations associated to this kind of algebras will contain matrices which
%are functions on  $\mathfrak{h}$. 
%For such a dynamical matrix $R(x)$ depending on $x$ taking values in $\textrm{End}(U_i \otimes U_j)$ we denote 
\egyenk\label{shift}
f(\lambda+\gamma h)=  f(\{\lambda_i+\gamma h_i\}) = \sum_{m \geq 0}\frac{\gamma ^m}{m!} \sum_{i_1, \dots ,i_m = 1}^n 
\frac {\partial^m f(\lambda)}{\partial \lambda_{i_1} \dots  
\partial \lambda_{i_m}}  h_{i_1} \dots h_{i_m} 
\egyenv
%The equations which we will write later 
%We will use the adjective ``dynamical'' as a synonim (?) for a ``function defined on $\mathfrak{h}$''. 
%If $U_1$ and $U_2$ are two finite dimensional semisimple $\mathfrak{h}$-modules then the ``shift'' of a function 
%\egyenk \label{shift}
%T: \mathfrak{h} & \to &\textrm{End} (U_1) \\
%x& \mapsto & T(x) \nonumber
%\egyenv 
%is defined as a function
%\egyenk
%\tilde{T}: \mathfrak{h} & \to &\textrm{End} (U_1 \otimes U_2)\\
%x& \mapsto & \tilde{T}(x)=
%T(x+\gamma h_2):= \nonumber \\
%&& = \sum_{n \geq 0}\frac{\gamma ^n}{n!} \sum_{i_1, \dots ,i_n = 1}^k 
%\frac {\partial_{i_1 \ldots i_n} T}{\partial \lambda_{i_1} \dots  
%\partial \lambda_{i_n}} \otimes h_{i_1} \dots h_{i_n} \nonumber
%\egyenv 
%A dynamical matrix will then be a function $R_{ij} : \mathfrak{h} \to \textrm{End} (V_i \otimes V_j)$.
which is a function on $\mathbb{C}^n$ taking values in $\mathsf{U}(\mathfrak{h})$.
%\texttt{needs to be clearer}

Armed with these definitions we propose the following dynamization of the algebra relations (\ref{QQA})
\egyenk
A_{12}(\lambda) T_1(\lambda) B_{12}(\lambda) T_2(\lambda+\gamma h_1) = T_2(\lambda) C_{12}(\lambda) 
T_1(\lambda+\gamma h_2) D_{12}(\lambda)\label{DQQA}
\egyenv
We require additional assumptions on the $R$-matrix
\egyenk \label{comm}
\left[h \otimes \mathbf{1}, B_{12} \right]=0,\quad 
\left[\mathbf{1} \otimes h,C_{12} \right]= 0,\quad
\left[h \otimes \mathbf{1}+\mathbf{1} \otimes h,D_{12} \right] =0 \quad (\forall h \in \mathfrak{h})
\egyenv
\egyenk
& B_{12} = C_{21}, \quad A_{21}=A_{12}^{-1}, \quad D_{21}=D_{12}^{-1}&  \label{flip}
\egyenv
%Assumptions (\ref{comm}) are motivated by the particular structure of the classical limit soon to be described.
Zero-weight conditions (\ref{comm}) will be presently seen to be consistent with the dynamical shifts in (\ref{DQQA}). 
We expect that different consistent choices of ``zero-weight conditions'' (\ref{comm}) will exist,
leading to different DQQA's, but we will not discuss it at this time.
Assumptions (\ref{flip}) are general self-consistency conditions for form-invariance of (\ref{DQQA}) under
exchange of labels $1$ and $2$.

Now the (sufficient) consistency conditions can be derived thanks to a change of point of view advocated in \cite{MFr}. 
Instead of looking
at $T$ as a matrix multiplied by, say, $A$ from left and $B$ from right one can think of $T$ as a bivector which is {\it acted upon} 
by $A$ and $B$. To put it another way: the triple matrix product $A \cdot T \cdot B$ can be viewed either as a sum 
$\ds\sum_{p,q}A^i_p T^p_q B^q_j$ -- where $T^p_q$ is viewed as a matrix element --, 
or as a sum $\ds\sum_{p,q} A^i_p \left( B^t \right) ^j_q T^{p,q}$ -- 
where $T^{p,q}$ is seen as a coordinate of a bivector of which each factor is multiplied separately.
Then we rewrite equation (\ref{DQQA}) as 
\egyenk
A_{12} T_{\bf{11'}}^{t_{1'}} B_{1'2} T_{\bf{22'}}^{t_{2'}}(\lambda+\gamma h_{1'}) = 
T_{\bf{22'}}^{t_{2'}} C_{12'} T_{\bf{11'}}^{t_{1'}}(\lambda+\gamma h_{2'}) D_{1'2'}\label{TT}
\egyenv
where bivector labels are indicated in bold. For the sake of simplicity, explicit dependence 
on dynamical parameters is omitted wherever possible. 
Using the commutation relations (\ref{comm}) (in fact here only the first two are needed) 
and the transposition on spaces $1'$ and $2'$ 
%yields:
%\egyenk
%A_{12} B_{1'2}^{t_{1'}} T_{\bf{11'}} T_{\bf{22'}}(\lambda+\gamma h_{1'}) = 
%D_{1'2'}^{t_{1'}t_{2'}}C_{12'}^{t_{2'}}T_{\bf{22'}} T_{\bf{11'}}(\lambda+\gamma h_{2'})
%\egyenv
%For the sake of simplicity, explicit dependence on dynamical parameters is omitted wherever possible.
equation (\ref{TT}) can  be recast into the form
\egyenk
\mathcal{R}_{11',22'} T_{\bf{11'}} T_{\bf{22'}}(\lambda+\gamma h_{1'}) = T_{\bf{22'}} T_{\bf{11'}}(\lambda+\gamma h_{2'})\label{ERv}
\egyenv
where $\mathcal{R}_{11',22'}$ is defined as
\egyenk
\mathcal{R}_{11',22'} = (C_{12'}^{t_{2'}})^{-1} (D_{1'2'}^{t_{1'}t_{2'}})^{-1}  A_{12} B_{1'2}^{t_{1'}}
\egyenv

The compatibility condition for the algebra generated by the elements of $T$ is derived as usual.
Starting from $$ T_{\bf{33'}} T_{\bf{22'}}(\lambda+\gamma h_{3'}) T_{\bf{11'}}(\lambda+\gamma h_{2'}+\gamma h_{3'})$$ 
one compares both ways (consistent by associativity) of obtaining 
$$ T_{\bf{11'}} T_{\bf{22'}}(\lambda+\gamma h_{1'}) T_{\bf{33'}}(\lambda+\gamma h_{1'}+\gamma h_{2'})$$ 
by permutation of the $T$-s using the
exchange relation (\ref{ERv}). 

\begin{L1}
A sufficient condition for the consistency of 3-space exchange of $T$-matrices with
\egyenk
\mathcal{R}_{11',22'} T_{\bf{11'}} T_{\bf{22'}}(\lambda+\gamma h_{1'}) = T_{\bf{22'}} T_{\bf{11'}}(\lambda+\gamma h_{2'})
\egyenv
is the following dynamical Yang-Baxter equation for $\mathcal R$:
\egyenk\label{gYB}
&& \mathcal{R}_{11',22'}(\lambda+\gamma h_{3'})  \mathcal{R}_{11',33'}(\lambda)  
\mathcal{R}_{22',33'}(\lambda+\gamma h_{1'}) \nonumber \\ &=& 
\mathcal{R}_{22',33'}(\lambda)  \mathcal{R}_{11',33'}(\lambda+\gamma h_{2'})  \mathcal{R}_{11',22'}(\lambda)
\egyenv
\end{L1}

Our goal is now to deduce from (\ref{gYB}) a set of consistent dynamical equations for the four components of 
the matrix $\mathcal{R}$. We will
illustrate this by explicitly describing the first step of the process in Appendix A. 
 
In the end we find that under 
assumptions (\ref{comm}) and (\ref{flip}) the nondynamical YB-equations obtained in \cite{MFr}
can be consistently dynamized as follows and this dynamization in turn assures that (\ref{gYB}) is satisfied.

\egyenk
 A_{12} \ A_{13} \ A_{23} \ &=& \ A_{23} \ A_{13} \ A_{12} \label{YBA}\\*
 D_{12}(\lambda+\gamma h_{3}) \ D_{13} \ D_{23}(\lambda+\gamma h_{1}) \ &=& \ D_{23} \ D_{13}(\lambda+\gamma h_2) 
\ D_{12} \label{YBD}\\*
 D_{12} \ B_{13} \ B_{23}(\lambda+\gamma h_1) \ &=& \ B_{23} \ B_{13}(\lambda+\gamma h_2) \ D_{12} \label{YBDB}\\*
 A_{12} \ C_{13} \ C_{23} \ &=& \ C_{23} \ C_{13} \ A_{12}(\lambda+\gamma h_3) \label{YBAC}
\egyenv

It can be checked that (\ref{YBA}) and (\ref{YBD})
are precisely the consistency conditions for the $BC$
algebras (\ref{YBDB}) and (\ref{YBAC}). For example starting with $$B_{14} \ B_{24}(\lambda+\gamma h_1) \ 
B_{34}(\lambda+\gamma h_1 + \gamma h_2)$$ and using
the exchange relation (\ref{YBDB}) one obtains $$B_{34} \ B_{24}(\lambda+\gamma h_3) \ 
B_{14}(\lambda+\gamma h_2 + \gamma h_3)$$ in two different
 ways. These two ways yield the same result whenever (\ref{YBD}) is satisfied. Note that (\ref{YBD}) together with 
the zero-weight condition (\ref{comm}) is the usual GNF equation. By contrast, $A$ obeys a non-dynamical Yang-Baxter equation a
lthough it also contains the dynamical variables.

%Where $x$ is a variable in the dual of a Cartan subalgebra $\mathcal H^{\star}$. If $T$ is a map $T :\mathcal {H^{\star}} 
%\to \textrm{End}(V)$ then its ``shift''  can be defined as 
%\egyenk
%T(x+\gamma h) = \sum_{n \geq 0}\frac{\gamma ^n}{n!} \sum_{i_1, \dots ,i_n = 1}^k \frac {\partial T}{\partial \lambda_{i_1} \dots  
%\lambda_{i_n}} \otimes h_{i_1} \dots h_{i_n}
%\egyenv

To summarize we now state 

\begin{T1} 
The exchange relations
\egyenk
A_{12}(\lambda) T_1(\lambda) B_{12}(\lambda) T_2(\lambda+\gamma h_1) = T_2(\lambda) C_{12}(\lambda) 
T_1(\lambda+\gamma h_2) D_{12}(\lambda) \nonumber
\egyenv
where
\egyenk \label{comm2}
\left[h \otimes \mathbf{1}, B_{12} \right]=0,\quad 
\left[\mathbf{1} \otimes h,C_{12} \right]= 0,\quad
\left[h \otimes \mathbf{1}+\mathbf{1} \otimes h,D_{12} \right] =0 \quad (\forall h \in \mathfrak{h})
\egyenv
\egyenk
& B_{12} = C_{21}, \quad A_{21}=A_{12}^{-1}, \quad D_{21}=D_{12}^{-1}&  \label{flip2}
\egyenv
  together with the relations 
 \egyenk
 A_{12} \ A_{13} \ A_{23} \ &=& \ A_{23} \ A_{13} \ A_{12} \label{thA}\\*
 D_{12}(\lambda+\gamma h_{3}) \ D_{13} \ D_{23}(\lambda+\gamma h_{1}) \ &=& \ D_{23} \ D_{13}(\lambda+\gamma h_2) 
\ D_{12} \label{thD}\\*
 D_{12} \ B_{13} \ B_{23}(\lambda+\gamma h_1) \ &=& \ B_{23} \ B_{13}(\lambda+\gamma h_2) \ D_{12} \label{thDB} \\*
 A_{12} \ C_{13} \ C_{23} \ &=& \ C_{23} \ C_{13} \ A_{12}(\lambda+\gamma h_3) \label{thAC}
\egyenv
yield an associative dynamical quadratic algebra.
\end{T1}
 
We now formulate two fusion structures on the quantum space.
\begin{T1}\label{Theorem2}
Let $T_{1q}$ be a representation of the algebra (\ref{DQQA}) on some Hilbert space
$H_q$. Let $L_{1q'}, R_{1q'}$ denote a representation on another Hilbert space $H_{q'}$ of the following 
set of exchange relations :
\egyenk
A_{12}\  L_1 L_2 &=& L_2 L_1 \ A_{12} \label{coprAB}\\
R_1 \  B_{12} \  L_2(\lambda+\gamma h_1) &=& L_2 \  B_{12} \  R_1 \nonumber\\
D_{12} \  R_1 R_2( \lambda+\gamma h_1) &=& R_2 R_1(\lambda+\gamma h_2) \  D_{12} \nonumber
\egyenv
then 
\egyenk T_{1,qq'}= L_{1q'} T_{1q} R_{1q'} \label{copr}\egyenv
 is a representation on $H_q \otimes H_{q'}$ of the algebra
(\ref{DQQA}).

Similarly, if $L_{1q'}, R_{1q'}$ is a representation on a Hilbert space $H_{q'}$ of the exchange relations
\egyenk
A_{12} \  L_1 L_2 &=& L_2 L_1 \  A_{12}(\lambda+\gamma h) \label{coprCD}\\
R_1 \  B_{12} \  L_2(\lambda+\gamma h_1) &=& L_2 \  B_{12}(\lambda+\gamma h) \  R_1 \nonumber\\
D_{12} \ (\lambda+\gamma h) R_1 R_2( \lambda+\gamma h_1) &=& R_2 R_1(\lambda+\gamma h_2) \  D_{12} .\nonumber
\egyenv
Then 
\egyenk T_{1,qq'}= L_{1q'} T_{1q}(\lambda+\gamma h_{q'}) R_{1q'} \label{copr2} \egyenv
yields a  representation of the same algebra (\ref{DQQA}) on the space $H_q \otimes H_{q'}$.
\end{T1}
It is assumed in (\ref{coprCD}) that the algebra of which $L_{1q'}$ and $R_{1q'}$ are representations
on $H_{q'}$ has an $\mathfrak{h}$-module structure, thereby making sense of the unindexed dynamical shift.

\vspace{3mm}

\textbf{Proof}: direct check of (\ref{DQQA}) by using the set of relations (\ref{coprAB}) or (\ref{coprCD}).

\vspace{3mm}
Strictly speaking we have here defined fusion procedures of represented $T$-matrices. 
In the sequel we will refer to these fusion structures simply as ``coproducts'' even though
we cannot prove yet that there is a universal bialgebra structure behind them.

A straightforward representation of the first $LR$-exchange algebra (\ref{coprAB}) is provided on $H_{q'}=V$ by taking
$L_{1q'} \equiv A_{12}$ and $R_{1q'} \equiv B_{12}$ or $L_{1q'} \equiv (A_{12}^{-1})^{t_2}$ 
and $R_{1q'} \equiv (B_{12}^{t_2})^{-1}$. The second $LR$ algebra, too, has a simple representation
in terms of the structure matrices. Namely, one can take $L_{1q'} \equiv C_{12}$ and $R_{1q'} \equiv D_{12}$
or $L_{1q'} \equiv (C_{12}^{t_2})^{-1}$ and $R_{1q'} \equiv (D_{12}^{-1})^{t_2}$. 
These representations by structure matrices are made possible by the fact that $A$ and $B$ (respectively $C$ and $D$) 
obey three out of the four consistency requirements: (\ref{thA}), (\ref{thDB}) and (\ref{thAC}) realizing (\ref{coprAB}) 
(respectively
(\ref{thD}), (\ref{thDB}) and (\ref{thAC}) realizing (\ref{coprCD})).
Let us remark here that both $LR$-algebras (\ref{coprAB}) and (\ref{coprCD}), therefore both coproducts, are identical 
in the nondynamical limit $\gamma \to 0$ 
to the single $T^{+},T^{-}$ algebra and its coproduct described in \cite{MFr}. This can be understood if one notices that
in the nondynamical case the consistency relations (\ref{YBA}) - (\ref{YBAC}) admit a particular symmetry: 
$A_{i3} \leftrightarrow C_{i3}$ and $D_{i3} \leftrightarrow B_{i3}$ for $i=1,2$.
 This is no longer true in the case $\gamma \neq 0$. 

\subsection{An example} 

A concrete realization of the algebra (\ref{DQQA}) is given \cite{russes} by the elliptic RS model \cite{RS}. 
For the sake of simplicity we only consider here its rational limit. Let us define the structure matrices as:
\egyenk
A(\lambda)=1+\sum_{i\neq j}\frac{\gamma}{\lambda_{ij}} \left(E_{ii}-E_{ij}\right) \otimes \left(E_{jj}-E_{ji}\right) \label{RSexA}\\
B(\lambda)=C(\lambda)^{\pi}=1+\sum_{i\neq j}\frac{\gamma}{\lambda_{ij}-\gamma} E_{jj} \otimes \left(E_{ii}-E_{ij}\right)\\
D(\lambda)=1-\sum_{i\neq j}\frac{\gamma}{\lambda_{ij}} E_{ii} \otimes E_{jj} + \sum_{i\neq j}\frac{\gamma}{\lambda_{ij}} E_{ij} 
\otimes E_{ji}\label{RSexD}
\egyenv
where $E_{ij}$ is the elementary matrix whose entries are $\left(E_{ij}\right)_{kl}=\delta_{ik}\delta_{jl}$ and 
$\lambda_{ij}=\lambda_i-\lambda_j$.
These matrices verify the consistency conditions (\ref{YBA})-(\ref{YBAC}). 
A scalar representation of the exchange algebra defined with these structure matrices is then provided by:
\egyenk
T(\lambda)=\sum_{ij} \frac{\prod_{a\neq i} (\lambda_{aj}+\tilde{\gamma})}{\prod_{a\neq j}\lambda_{aj}} E_{ij} \otimes \mathbf{1}
\label{RSexL}
\egyenv

Taking this as a starting point, one can now use the coproducts described above to construct other, 
higher dimensional nonabelian 
representations of the algebra defined by 
(\ref{RSexA}) - (\ref{RSexD}) which should provide us with a suitable algebraic framework to define and study 
spin generalizations of the RS-model.  
Theorem \ref{Theorem2} indeed provides us with the dynamical version of the construction
of a monodromy matrix for a spin chain model by successive products of $R$-matrices, using the
coproduct structure of the quantum group.

\section{The classical limit}

For classical integrable systems the starting point is the following quadratic Poisson-bracket algebra \cite{MFr, Mail2}
\egyenk
\{l_1,l_2\}=a_{12} l_1 l_2+l_1 b_{12} l_2-l_2 c_{12} l_1 -l_1 l_2 d_{12}
\egyenv
where the Lax-matrix $l$ is a function on the phase space taking values in $\textrm{End}(V)$, $V$ being a finite dimensional 
vector space. The matrices $a$, $b$, $c$, $d$ that define the quadratic algebra 
are elements of $\textrm{End}(V \otimes V)$. We say that the algebra is dynamic if these matrices actually 
depend on the phase space variables. 

In order to ensure the antisymmetry of the Poisson-bracket we impose the following constraints on the structure matrices:
\egyenk
a+a^{\pi} =\alpha \mathcal{C}, \qquad d+d^{\pi}=\alpha \mathcal{C}, \qquad b^{\pi}=c \qquad \ (\alpha \in \mathbb{C})
\egyenv
where, as usual, $\pi$ denotes the permutation in $\textrm{End}(V \otimes V)$, and  
$\mathcal{C}$ is the Casimir-operator, i.e. for the $\mathfrak{gl}_n$ case $\mathcal{C}=\sum_{i,j} E_{ij} \otimes E_{ji}$.
In other words, we are allowed to modify $a$ and $d$ by adding the same multiple of $\mathcal{C}$ to both of them. The Poisson-bracket
will not change, since $\left[\mathcal{C},l_1 l_2\right]=0$. 
These conditions on $a$ and $d$ are slightly more relaxed than usual: the reason for this will become clear when
we consider the RS model. 

A well-behaved Poisson-bracket should also verify the Jacobi identity. This is equivalent to demanding that the following general
identity holds:
\egyenk
&\left( \left[a_{12},a_{13}\right]  +  \left[a_{12},a_{23}\right]  +  \left[a_{13},a_{23}\right]\right)l_1 l_2 l_3 - \nonumber \\*
&-l_1 l_2 l_3 \left( \left[d_{12},d_{13}\right]  +  \left[d_{12},d_{23}\right]  +  \left[d_{13},d_{23}\right] \right) + \nonumber \\*
&+l_1 l_2 \left( \left[d_{12},b_{13}\right]+\left[d_{12},d_{23}\right]+\left[b_{13},b_{23}\right] \right) l_3 + \textrm{circ. perm.}
\nonumber \\*
&-l_3 \left( \left[a_{12},c_{13}\right]+\left[a_{12},c_{23}\right]+\left[c_{13},c_{23}\right] \right) l_1 l_2 - \textrm{circ. perm.}
\nonumber \\*
&-l_1 l_2 \{d_{12},l_3\} - l_2 l_3 \{d_{23},l_1\} - l_3 l_1\{d_{31},l_2\} \nonumber \\*
&+\{a_{12},l_3\} l_1 l_2 + \{a_{23},l_1\} l_2 l_3 + \{a_{31},l_2\} l_1 l_3 \nonumber \\* 
&-l_2 \{c_{12},l_3\} l_1 - l_3 \{c_{23},l_1\} l_2 - l_1 \{c_{31},l_2\} l_3 \nonumber \\*
&+l_1 \{b_{12},l_3 \} l_2 + l_2 \{b_{23},l_1\} l_3+ l_3 \{b_{31},l_2\} l_1 = 0 \label{PB}
\egyenv
where the last four lines appear because of the dynamical nature of the structure matrices.

\subsection{An example: the hyperbolic Ruijsenaars-Schneider model}

Due to the appearance of dynamical terms of the generic form $\{a,L\}$ in the Jacobi identity, 
it is not clear how to characterize general algebraic structures in (\ref{PB}). 
However, again in the concrete example of the  Ruijsenaars-Schneider model, the particular form
of the occuring matrices enables us to proceed a step further, define a fully algebraic classical YB formulation
 and eventually connect it with our dynamical quadratic quantum algebras. 
Let us consider the Lax-matrix structure in the RS $\textrm{A}_n$ 
case which reads as follows \cite{Suris,AR}:
\egyenk
&l=\sum_{i,j} l_{ij} E_{ij}& \\*
&l_{ij}=c(q_i-q_j) e^{-p_j} f_j& 
\egyenv
where $c$ and $f_j$ are functions of the position variables.
The Poisson-bracket on the phase space is given by
$\{p_i,q_j\}=\delta_{ij}$.
 
The quadratic structure coefficients read \cite{Suris}:
\egyenk
&&a = -u -s +s^{\pi} +w-\mathcal{C}, \quad b=-s^{\pi}-w, \nonumber \\
&&c= -s+w, \quad d=-u-w-\mathcal{C} \label{RSr}
\egyenv
For the hyperbolic model the matrices $u,s,w$ take the form:
\egyenk
&&u=-\sum_{i\neq j}  \coth (q_i-q_j)\ E_{ij} \otimes E_{ji}\ , \quad s=\sum_{i\neq j}\frac{1}{\sinh (q_i-q_j)}\  E_{ij}
\otimes E_{jj} \nonumber \\*
&&w=\sum_{i\neq j}  \coth (q_i-q_j) \ E_{ii}\otimes E_{jj}
\egyenv
Now
using the fact that the $a$, $b$, $c$, $d$ matrices depend only on the position variables and 
that $L_{ij}$ depends on $p$ as $e^{-p_j}$, the Poisson-brackets  in the last four lines of (\ref{PB}) can be 
written as $L$ multiplying a certain sum {\it from the left}: 
\egyenk
\{M_{12},l_3\} = l_3 \sum_k E_{kk}^{(3)} \partial_k M_{12}
\egyenv
where $M$ stands for any matrix depending only on the position variables. As a result we can rewrite the Jacobi-identity 
in a purely algebraic form as follows.
\egyenk
&\Big( \left[a_{12},a_{13}\right]  +  \left[a_{12},a_{23} \right] + \left[a_{13},a_{23}\right] \Big)\ l_1 l_2 l_3 - \nonumber \\*
&-l_1 l_2 l_3 \ \Big( \left[d_{12},d_{13}\right]  +  \left[d_{12},d_{23}\right]  +  \left[d_{13},d_{23}\right] + \nonumber \\*
& + \sum_k h^{(1)}_k \partial_k d_{23} - \sum_k h^{(2)}_k \partial_k d_{13} + \sum_k h^{(3)}_k \partial_k d_{12} \Big)   \nonumber \\*
& +l_1 l_2\ \Big( \left[d_{12},b_{13}\right]+\left[d_{12},d_{23}\right]+\left[b_{13},b_{23}\right] + \nonumber \\*
& + \sum_k  h^{(1)}_k \partial_k  b_{23} - \sum_k h^{(2)}_k \partial_k b_{13} \Big)\ l_3+ \textrm{circ. perm.} \label{PB2} \\*
& - l_3\ \Big( \left[a_{12},c_{13}\right]+\left[a_{12},c_{23}\right]+\left[c_{13},c_{23}\right] -\sum_k h^{(3)}_k \partial_k a_{12} 
\Big)\ l_1 l_2- \textrm{circ.perm.} = 0\nonumber
\egyenv

Hence we may now state:

\begin{T1}
A set of sufficient conditions on the $r$-matrix which ensures that the Jacobi-identity  holds is
\egyenk
&\left[a_{12},a_{13}\right]  +  \left[a_{12},a_{23}\right]  +  \left[a_{13},a_{23}\right] = 0 \label{cYBA}\\*
&\left[d_{12},d_{13}\right]  +  \left[d_{12},d_{23}\right]  +  \left[d_{13},d_{23}\right]  +  \sum_k h^{(1)}_k \partial_k d_{23} 
 -  \sum_k h^{(2)}_k \partial_k d_{13}+ \label{cYBD}\\*
&+  \sum_k h^{(3)}_k \partial_k d_{12} = 0 \nonumber \\* 
&\left[d_{12},b_{13}\right]+\left[d_{12},b_{23}\right]+\left[b_{13},b_{23}\right] +\sum_k  h^{(1)}_k \partial_k  b_{23} - 
\sum_k h^{(2)}_k \partial_k b_{13}=0 \label{cYBDB}\\*
&\left[a_{12},c_{13}\right]+\left[a_{12},c_{23}\right]+\left[c_{13},c_{23}\right]-\sum_k h^{(3)}_k \partial_k a_{12}=0 \label{cYBAC}
\egyenv
\end{T1}

We are now able to establish a link between  these equations and the quantum algebra presented in (\ref{DQQA}). 
Indeed, if we assume the existence of a classical limit for $A,B,C,D$ in (\ref{YBA})-(\ref{YBAC}) as 
$A=\mathbf{1} + \hbar \  a +\textrm{O}(\hbar^2)$ ... 
%If we suppose that the $r$-matrix in (\ref{RSr}) is a quasiclassical limit of a quantum $R$-matrix, i.e. $A=\mathbf{1} + \gamma \  a 
%+\textrm{O}(\gamma^2)$, etc.,  
we can expand the quantum YB-equations (\ref{YBA})-(\ref{YBAC}) in powers of $\hbar$. 
The equations (\ref{cYBA})-({\ref{cYBAC}) will appear  as the first nontrivial term  
(of order $\hbar ^2$) in this expansion.

An example of solution to these classical quadratic YB-equations (\ref{cYBA}) - (\ref{cYBAC}) is again provided by \cite{russes}:
\begin{P1}
 The $r$-matrix of the RS-model defined in (\ref{RSr}) verifies
these equations. 
\end{P1}

This is easily established by direct computations, and does require the additional Casimir terms in $a$ and $d$.
This provides us in addition with a complete algebraic interpretation of the YB-equation for the non-antisymmetric
$r$-matrix of the scalar CM-model obtained as $r=a-c=d-b$ (up to the extra Casimir terms). 
%This interpretation was lacking until now.
The full antisymmetric part $d$ and the non-antisymmetric part $b$ must be treated as separate objects obeying
(\ref{cYBD}) and (\ref{cYBDB}). This of course explains why such $r$-matrices are absent from the classification in 
\cite{Schiffmann2} where only solutions to (\ref{cYBD}) are considered.
%Note also that the $r$-matrix in (\ref{RSr})
%satisfies the commutation relations (\ref{comm}); the procedure is therefore completely consistent.

\section{Conclusion}

We have proposed a dynamical extension for a general quadratic algebra;
we have explicited the consistency conditions as a set of dynamical YB-type equations
generalizing the set given in \cite{MFr}, and we have constructed two independent coproduct structures for them.
%At this time we have no representation of the
%quantum equations but we know at least one realization of its semi-classical limit, e.g. the
%$r$-matrix structure for the $A_n$ Ruijsenaars-Schneider Lax representation.

The next steps are clearly defined: we will look for new explicit quantum solutions of the set (\ref{YBA})-(\ref{YBAC})
 by combining the initial representation (\ref{RSexA})-(\ref{RSexL}) with the coproduct (\ref{copr}) which should lead
us to ``spin-RS''-like Lax matrices for which this structure would thus provide a suitable algebraic framework;
and we will also look for other consistent dynamical extensions of the quadratic algebras.
Indeed one already knows at least two such structures: the quantum dynamical Gervais-Neveu-Felder algebras \cite{Fe,GN}
, where $B=C=\mathbf{1}$ and $A=D$ with $h_1+h_2$ zero-weight condition;
and a suggested dynamical version of the reflection algebras, where $A=D^{\pi}$, $B=C^{\pi}$ and all objects derive 
from a single spectral parameter dependent $R$-matrix \cite{Ra}. It would be very significant to understand the general scheme, 
and we hope to report on this soon.

Another issue would be the understanding of (\ref{YBA}) - (\ref{YBAC}) as defining relations for some (quasi-Hopf?) bialgebra
, generalizing the construction of \cite{DKM} for the nondynamical case by suitably incorporating the ``coproducts'' (\ref{copr}),
 (\ref{copr2}).

\vspace{5mm}
 {\bf Acknowledgments}

We wish to thank J.~M.~Maillet for his explanations on the bivector method described in \cite{MFr}; P.~P.~Kulish 
and M.~Jimbo for enlightening suggestions; E.~Ragoucy for pointing out 
to us ref. \cite{Ra}, A. Doikou for clarifying comments and showing us refs. \cite{McGu, Ber,  BZinn}.
\setcounter{section}{0}
\appendix{}

We start by writing down the left hand side of equation (\ref{gYB}) in its full length.

\egyenk 
&(C_{12'}^{t_{2'}})^{-1}(h_{3'}) (D_{1'2'}^{t_{1'}t_{2'}})^{-1}(h_{3'})  A_{12}(h_{3'}) B_{1'2}^{t_{1'}}(h_{3'}) \times\\ 
&(C_{13'}^{t_{3'}})^{-1} (D_{1'3'}^{t_{1'}t_{3'}})^{-1}  A_{13} B_{1'3}^{t_{1'}}  \times \nonumber\\ 
&(C_{23'}^{t_{3'}})^{-1}(h_{1'}) (D_{2'3'}^{t_{2'}t_{3'}})^{-1}(h_{1'})  A_{23}(h_{1'}) B_{2'3}^{t_{1'}}(h_{1'})  \nonumber
\egyenv

We then pick $B_{1'2}^{t_{1'}}(h_{3'})$ and push it through the rest of the product as far as the commutation relations 
(\ref{comm}) make it possible: in our case $( D_{1'3'}^{t_{1'}t_{3'}} )^{-1}$ blocks the way. Now we select the matrix 
acting on spaces that allow to form a YB-type
equation with $B_{1'2}^{t_{1'}}(h_{3'})$ and $( D_{1'3'}^{t_{1'}t_{3'}} )^{-1}$. 
In this case it is $(C_{23'}^{t_{3'}}(h_{1'}))^{-1}$. We indeed can push it through
the matrices separating from $( D_{1'3'}^{t_{1'}t_{3'}} )^{-1}$ thanks again to the commutation relations (\ref{comm}). 
We now have to fix 
a suitable exchange relation for $B$, $D$ and $C$. A consistent choice is:
\egyenk\label{YB1}
B_{1'2}^{t_{1'}}(h_{3'}) \ ( D_{1'3'}^{t_{1'}t_{3'}} )^{-1} \ (C_{23'}^{t_{3'}}(h_{1'}))^{-1} = 
(C_{23'}^{t_{3'}})^{-1}\ (D_{1'3'}^{t_{1'}t_{3'}})^{-1}\  B_{1'2}^{t_{1'}} 
\egyenv
which yields after some rearranging, using (\ref{flip}), and total transposition:
\egyenk
B_{13}^{t_{3}}(h_{2})\ B_{23}^{t_{3}} \ D_{12} = 
D_{12} \ B_{23}^{t_{3}}(h_{1})\ B_{13}^{t_{3}}
\egyenv
Thanks to the commutation relations (\ref{comm}), it is possible to transpose this equation on space $3$.
%Indeed, let us write down a general term of the sum (\ref{shift}):
%\egyenk
%\partial^I B_{13}^{t_{3}}\ h_{2}^I \  B_{23}^{t_{3}} \  D_{12} = \partial^I B_{13}^{t_{3}}\ B_{23}^{t_{3}} h_{2}^I \ D_{12} \nonumber
%\egyenv
%where $\partial^I$ and $h_{2}^I$ are shorthand notations\footnote{$h_{2}^I$ should be noted $h_{I2}$, but we chose 
%to raise $I$ in order to avoid
%confusion between indices that label spaces and elements of a basis.} for $\partial_{\lambda_{i_1} \ldots \lambda_{i_n}}$ a
%nd $\mathbf{1} 
%\otimes h_{i_1} \cdots h_{i_n} \otimes \mathbf{1}$.
%After transposition on space $3$ this becomes
%\egyenk
%B_{23}\ \partial^I B_{13} h_{2}^I \ D_{12}\nonumber
%\egyenv
%Which is the general term of 
%\egyenk
%B_{23} \ B_{13}(h_2) D_{12}
%\egyenv
%Then transposing on space $1$  yields
%\egyenk
%h_{2}^I \ B_{23}\ \partial^I B_{13}\ D_{12} = B_{23}\ \partial^I  B_{13}\ h_{2}^I\ D_{12} \nonumber
%\egyenv
%where we again require the zero-weight property (\ref{comm}). 
The right hand side is to be treated similarly. 
Eventually we find that requiring (\ref{YB1}) amounts to imposing the following dynamical YB-equation on $D$
and $B$:
\egyenk
D_{12} B_{13} B_{23}(\lambda+\gamma h_1) = B_{23} B_{13}(\lambda+\gamma h_2) D_{12}
\egyenv
The remaining YB-equations are obtained by repeating the above described process.

\end{document}